\documentclass[12pt,a4paper,oneside]{article}
\usepackage{amsmath,amstext,amssymb,amsopn,amsthm,color}
\theoremstyle{plain}
\newtheorem{thm}{Theorem}[section]

\newtheorem{lem}[thm]{Lemma}

\newtheorem{cor}[thm]{Corollary}

\theoremstyle{definition}

\theoremstyle{remark}
\newtheorem*{rem*}{Remark}

\newcommand{\R}{\mathbb{R}}
\newcommand{\N}{\mathbb{N}}

\newcommand{\C}{\mathbb{C}}
\renewcommand{\H}{\mathbb{H}}

\renewcommand{\leq}{\leqslant}
\renewcommand{\geq}{\geqslant}

\newcommand{\pref}[1]{(\ref{#1})}

\def\o{\over}
\def\({\left(}
\def\){\right)}
\def\[{\left[}
\def\]{\right]}
\def\<{\langle}
\def\>{\rangle}

\title {Hitting distributions of geometric Brownian motion
\footnotetext{2000 MS Classification:
    Primary 60J65;
    Secondary 60J60. {\it Key words and phrases}: hyperbolic spaces, stable
    processes, Poisson kernel.  Research partially supported by KBN
    grant 1 P03A 020 28 and RTN Harmonic Analysis and Related Problems
    contract HPRN-CT-2001-00273-HARP}}
\author{ T. Byczkowski and M. Ryznar\\ Institute of Mathematics, Wroc\l{}aw
    University of Technology, Poland}
\date{}

\begin{document}
\maketitle

\begin{abstract}
    Let $\tau$ be the first hitting time of the point $1$ by the geometric Brownian motion
   $X(t)= x \exp(B(t)-2\mu t)$ with drift $\mu \geq 0$ starting from $x>1$.
   Here $B(t)$ is the Brownian motion starting from $0$ with $E^0 B^2(t) = 2t$.
   We provide an integral formula for the density function of the stopped exponential
  functional $A(\tau)=\int_0^\tau X^2(t) dt$ and determine
  its asymptotic behaviour at infinity. Although we basically rely on methods developed
  in \cite{BGS}, the present paper also covers the case of arbitrary drifts
   $\mu \geq 0$ and
  provides a significant unification and extension of results
  of the above-mentioned paper.  As a corollary we provide an integral formula and
  give asymptotic behaviour at infinity of the Poisson kernel for half-spaces
  for Brownian motion with drift in real hyperbolic spaces of arbitrary dimension.

\end{abstract}

\newpage
\section{Introduction}

Let $B(t)$ be the standard Brownian motion starting from $0$ and normalized
such that $E^0 B(t)^2 = 2t$. Consider the following linear SDE:
$$dX(t)= X(t) dB(t) - (2\mu -1) X(t) dt\/, \quad X(0)=x>0\/, \quad \mu \geq 0\/.$$
The strong unique non-exploding solution is given by
\begin{equation*}
X(t)=x \exp(B(t) - 2\mu t).
\end{equation*}
The process $\{X(t); t\geq 0\}$ is called a {\it geometric Brownian motion}
or {\it exponential Brownian motion} and  along with
 the additive functional from the process $X(t)$,
\begin{equation*}
A_x(t)= \int_0^t X^2(t) dt = \int_0^t x^2 \exp 2(B(s)-2\mu s)\/ ds
\end{equation*}
is of primary interest in mathematical
finance and insurance theory (see, e.g., \cite{D}, \cite{GY} or \cite{Y2}).
 Also there is a connection of the above functional with the Brownian motion
 in hyperbolic halfspaces (see eg. \cite {Y2}, \cite {AG}, \cite {BTF}, \cite {BTFY}
  and Section 5).

The distribution of  $A_x(t)$ for fixed $0\le t\le \infty$ has been a subject of study
 in substantial number of papers (see e.g. \cite {B}, \cite {Y2}, \cite {AG}, \cite {M}).
In this paper we investigate the properties of the
density function $q_\mu$ of the (stopped) additive functional $A_x(\tau)$
where $\tau$ is the first hitting time of the point $1$ by the process $X(t)$
(starting from $x>1$). From the strong Markov property it easily follows that
the distribution of $A_x(\tau)$ is closely related to the distribution of $A_x(\infty)$.
  The latter one is astonishingly simple;
it is identical
  with
the distribution of the random variable $x^2/4Z$, where $Z$ is a
$\Gamma(\mu,1)$-distributed random variable (with the density
$\Gamma(\mu)^{-1}\/ u^{\mu-1} e^{-u}$).
This fundamental result is due to D. Dufresne \cite{D} (see also \cite{Y1},
 \cite{Y2}) and is of primary importance here. Knowing the  Laplace transform of
  $A_x(\infty)$ one can easily  derive  (via the strong Markov property) the
form of the Laplace transform of the random variable $A(\tau)$ as a suitable ratio
 of Bessel functions. For the reader's convenience we present this argument
  in Preliminaries.

On the other hand the distributions of $A_x(\tau)$  and $A_x(\infty)$ are
 closely related to hitting  times of Bessel processes.
  The main fact here is Lamperti's representation which reads  that there
   exists a Bessel process $R^{(-\mu)}$ with index $-\mu$, starting from $x>1$,
such that  the process $X(t)$ admits the following representation (see \cite{L}
 and  Ex. 11.1.28 in \cite {RY}):
\begin{equation*}
X(t) = R^{(-\mu)}\left(A_x(t)\right), \quad t\ge 0.
\end{equation*}
We refer the reader to  \cite {RY} for the account on  Bessel processes (see also
\cite {GJY}). From the Lamperti representation   it follows immediately that
 $A_x(\tau)$ and $A_x(\infty)$ can be regarded as  hitting times of $1$ and $0$,
  respectively,  of the Bessel process $R^{(-\mu)}$ starting from $x>1$. Also it is
   possible to relate $A_x(\tau)$ and  $A_x(\infty)$ to last exit times for
    appropriate Bessel processes. There is a vast existing literature on that subject.
     Such hitting or last exit times were studied in papers by Getoor \cite{G},
      Getoor-Sharpe \cite{GS}, Kent \cite{K} and Pitman-Yor \cite{PY}. In Getoor
       \cite{G}, Kent \cite{K} and Getoor-Sharp \cite{GS} the Laplace transform of
        the distribution of $A(\tau)$ is derived as a ratio of Bessel functions
(see Preliminaries). For an exhausting discussion on that subject we refer
to \cite {GJY}.

Our main focus in this paper is to provide an integral formula for the density of
 $A_x(\tau)$ (see Theorem \ref {rep}). In the case $\mu =1/2$ this density is well
  known to be $1/2$-stable subordinator. Kent in \cite K writes that it is possible
   to obtain the explicit form of the density for $\mu=3/2$ but he does not provide
    any formula nor  details on that.

At this point let us mention that various ratios of Bessel functions were proved to be
 completely monotonic functions hence they are Laplace transforms  of probability
  distributions (see Ismail \cite{I1},\cite{I2}, Ismail-Kelker \cite{IK},
   Kent \cite{K}). For a survey on this theme see Pitman and Yor \cite{PY} and
    also \cite {GY}. Regarding the case which is considered in our paper
     Ismail and Kelker \cite{IK} showed by purely analytical methods that it is  an
      infinitely divisible distribution.

The main purpose of the paper is to obtain a suitable representation for the density
 function
$q_\mu$ of the functional $A_x(\tau)$ along with its asymptotic properties, for
 arbitrary
drift terms $\mu \geq 0$. Then we apply these results to derive an integral
 representation of the Poisson kernel for subspaces on real hyperbolic spaces,
  for hyperbolic Brownian motion with arbitrary drift, which extends and simplifies
   the results and proofs from \cite{BGS}.

The paper is organized as follows. In Preliminaries we collect basic information
needed in the sequel.

 In Section 3 we obtain a suitable representation of the density
 of the functional $A_x(\tau)$. For this purpose we apply an integral
  representation given in \cite{BGS} for ratio of Bessel functions.

In Section 4 we exhibit the exact asymptotics of the density
of $A_x(\tau)$ at infinity for all drifts $\mu\geq 0$. Again, we essentially follow
the idea of \cite{BGS}.
However, applying more direct probabilistic arguments, we are able to simplify
our presentation considerably.

 In Section 5 we show how to apply results obtained in preceding
sections to obtain representation and asymptotic properties of Poisson kernel
 of subspaces
on real hyperbolic spaces of arbitrary dimension, for Brownian motion with drift.
Some of these asymptotic properties were studied in \cite{BTF} in the case
 of dimension 2, where it was shown that the distribution of  Poisson kernel belongs
  to the stable  domain of attraction. Our asymptotic results may be viewed as
   the extension of those obtained in \cite{BTF}. For a related result see
    also \cite{BTFY}.

 \section{Preliminaries}

 Let $0<a<x$ and let $\tau_a$ be the first hitting time of the point $a$ by the
 geometric Brownian motion with drift $\mu \geq 0$ starting at $x$:
 \begin{equation*}
 \tau_a = \inf\{t>0; x\exp(B(t)-2\mu t)=a\}.
  \end{equation*}
The fact that $\tau <\infty$ a.e. follows from
the property $\inf_{t>0} B(t) = - \infty$.

Further, define
\begin{equation*}
A_x(t)=\int_0^t \exp 2(B(t)-2\mu t)\/dt.
\end{equation*}
By the strong Markov property of Brownian motion we obtain

 \textbf{Basic relationship} (for $\mu >0$).
 \newline Observe that $x^2\exp 2(W(\tau_a)-2\mu \tau_a)=a^2$, hence
  \begin{eqnarray}
  A_x(\infty)&=&x^2\int_0^{\tau_a}\exp 2(W(s)-2\mu s)\/ds +
  x^2\int_{\tau_a}^\infty\exp 2(W(s)-2\mu s)\/ds \nonumber \\
  &=&  A_{x}(\tau_a)+ x^2\exp 2(W(\tau_a)-2\mu \tau_a)
  \int_0^\infty \exp 2(W(s+\tau_a)-W(\tau_a)-2\mu s)\/ds\nonumber\\
 & =&A_{x}(\tau_a)+  A'_a(\infty)\label{basic},
  \end{eqnarray}
  where $A'_a(\infty)$ is a copy of $A_a(\infty)$, independent from $A_{x}(\tau_a)$.

Dufresne \cite{D} (see also   Getoor \cite{G}, Kent \cite{K}, Getoor-Sharpe \cite{GS},
 where the result is given in the context of Bessel processes) showed that the density
  and the Laplace transform of $A_1(\infty)$ have the following form:
   \begin{equation}\label{subord}
    h_{\mu}(t)={2^{-2\mu} \o \Gamma(\mu)}  {{e^{-1/4t}} \o t^{1+\mu}},
    \end{equation}
  \begin{equation}\label{LAinfty}
  {\cal{L}}\{ A_1(\infty)\}(r^2)={\cal{L}}\{1/4Z\}(r^2)=
  2{(r/2)^\mu \o \Gamma(\mu)} K_\mu(r)\/, \quad \mu >0\/.
  \end{equation}
 From this and elementary properties of Laplace transform one immediately  obtains
  \begin{equation}\label{LAtau}
    {\cal{L}}\{ A_{x}(\tau_a)\}(r^2)=
   \({x \o a}\)^\mu {K_\mu(xr) \o K_\mu(ar)}\/,
    \end{equation}
where $K_\mu$ is the modified Bessel function of second type with index $\mu$.
The above Laplace transform appears in many papers (see eg.  Kent \cite{K},
 Getoor-Sharpe \cite{GS}, Ismail-Kelker \cite{IK}) and here we followed a
  derivation from \cite {BTF}. As mentioned in the Introduction we present it
   for the reader's convenience.

From the continuity of both sides of \pref{LAtau} with respect to
$\mu \geq 0$ we obtain that the above formula remains valid also
for $\mu=0$.

 Let us remark that for $\mu=0$ the right-hand side of \pref{LAtau} gives also
 the classical
 formula (\cite{S}) for the Laplace transform of Brownian motion
  hitting time $T_a$ of the centered circle
 with radius $a$ in $\R^2$ from the point $y \in \R^2$ such that $|y|=x>a$.

 By the formula \pref{LAtau} follows directly that for $t>0$ we have
 \begin{equation*}
 A_{ta,tx}(\tau_{ta}) \stackrel{d}= t^2 A_{a,x}(\tau_a) ,
  \end{equation*}
 where the symbol $\stackrel{d}=$ denotes the equality of distribution.

 Therefore, from now on we may and do assume
 that $a=1$ and $x>1$ is fixed. We write $\tau$ instead of $\tau_1$ and
 $A(\tau)$ instead of  $A_{1,x}(\tau_1)$.

  We conclude this section with the following technical lemma

 \begin{lem}\label{Pkernel}

   \begin{equation} \label{kernel}
   \int_0^\infty e^{-y^2/4t} h_s(t/\lambda^2) dt =
   {\lambda^{2s+2} \o (\lambda^2 + y^2)^{s}}\/.
   \end{equation}
  \end{lem}
  \begin{proof} Indeed, we obtain
  \begin{eqnarray*}
  &{}& \int_0^\infty e^{-y^2/4t} e^{-\lambda^2/4t} {\lambda^{2s+2}dt \o t^{1+s}}
     ={\lambda^{2s+2} \o (\lambda^2 + y^2)^s} \int_0^\infty e^{-(y^2+\lambda^2)/4t}
       {dt/(y^2+\lambda^2) \o (t/(y^2+\lambda^2))^{1+s}} \\
  &=& {2^{2s} \lambda^{2s+2} \o (\lambda^2+ y^2)^{s}} \int_0^\infty \Gamma(s) h_s(u) du
          = { 2^{2s}\Gamma(s) \lambda^{2s+2} \o (\lambda^2+ y^2)^{s}}\/.
   \end{eqnarray*}
  \end{proof}

Throughout the paper we use the following convention: by $c,\ C$ we always denote
 nonnegative constants which may depend on other constant parameters only.
  The value of $c$ or $C$ may change from line to line in a chain of estimates.

The notion $p(u)\approx q(u),\ u \to u_0$ means that the ratio $q/p\to 1$ when
$u \to u_0$.
\section{Representation of density of $A(\tau)$}\label{pkhf}

In this section we give a  representation formula for the density function of
functional $A(\tau)$, for arbitrary $\mu \geq 0$.

 From now on we use the following  notation, partially introduced in the preceding Section:
 \begin{equation*}
 x>1,\ \ a=1, \ \
 \lambda=x-1\/.
 \end{equation*}
and we denote by $q_\mu$ the density function of the functional
$A(\tau)$.

We begin with stating a more general version of a lemma, taken from \cite{BGS}.
The proof is identical and is omitted.

\begin{lem}\label{dzielenie}
  Let $\mu \geq0$. Suppose that
  \begin{equation*}
    Q(z) = z-(\mu^2-1/4) {\lambda \o 2x},\quad z\in \C.
  \end{equation*}
 Define  $F_\lambda(z)$ by the following formula:
  \begin{equation} \label{Flambda0}
    \lambda F_\lambda(z) =
    {z e^{\lambda z}x^{\mu} K_{\mu}(xz) - x^{\mu -{1\o2}}
    Q(z) K_{\mu}(z) \o    K_{\mu}(z)}.
  \end{equation}
  Then
  \begin{equation}\label{zanik}
    F_\lambda(z) = O(z^{-1}),\quad  z\to \infty\/.
     \end{equation}
 and there exists a function $w_\lambda$ such that
  \begin{equation}\label{laplace0}
      F_\lambda(z) = \int_0^\infty e^{-zv} w_\lambda(v) dv\/.
      \end{equation}
Moreover,
   \begin{equation}\label{laplace01}
   x^{\mu-1/2} (\mu^2-1/4)/2x=\int_0^\infty w_\lambda(v) dv\/,
   \end{equation}
   and, for $\mu>1/2$,
 \begin{equation}\label{laplace02}
  2 x^{\mu-1/2} =\int_0^\infty \kappa w_\lambda(v) dv,
   \end{equation}
 where $\kappa = (\lambda+v)^2 - \lambda^2 = v(2\lambda+v).$

   For $ \mu=1/2$ we have  $F_\lambda(z)\equiv 0$.
 \end{lem}
An explicit formula for the function $w_\lambda$ will be provided
in the sequel.

The following formula is crucial for our considerations:

 \begin{lem} \label{costrans}
Let $\varphi : [0,\infty) \rightarrow \R^{+}$ belong to $L^1(0,\infty)$
 such that its Laplace transform $\hat{\varphi}$ has the following property
 \begin{equation*} 
 t^{-1/2}  \hat{\varphi}(t) \in L^1(0,\infty).
 \end{equation*}
 Then we have
 \begin{equation} \label{expintegr1}
  \int_0^\infty \hat{\varphi}(r^2) \cos(ry) dr =
 \sqrt{\pi}/2 \int_0^\infty  e^{-y^2/4t} \varphi(t) { dt \o \sqrt{t}}.
  \end{equation}
 \end{lem}
 \begin{proof}
 Let $A$ be a random variable with absolutely continuous distribution with the
 density function $\varphi$ and let $W(t)$ be Brownian Motion starting from $0$
 (such that $EW^2(t)=2t$), independent from $A$. It is easy to see that the value
 of the Laplace transform of $A$ at the point $r^2$, that is $\hat{\varphi}(r^2)$
 is equal to the Fourier (or cosine transform) of $W(A)$ at the point $r$. Observe
 that our assumption assures that this Fourier transform belongs to $ L^1(0,\infty)$.
 Thus, the left-hand side of \pref{expintegr1} is the inversion formula applied
 for the Fourier transform and gives the density of $W(A)$. The right-hand side
 results from the direct computation of this density, taking into account the
 independence of $W$ and $A$ and the particular form of the (gaussian) density
 of $W(t)$.
  \end{proof}

We are ready to state our representation formula.
\begin{thm}\label{rep}
  For $\mu \geq 0$ we have
  \begin{equation}\label{rep1}
    q_{\mu}(t) = \lambda {e^{-\lambda^2/4t} \o \sqrt{\pi t} }
   \(x^{\mu-1/2}/2t +
    \int_0^\infty \(e^{-\kappa/4t} - 1\) w_\lambda(v)dv \) \/.
    \end{equation}

For $\mu>1/2$ we have
\begin{equation}\label{rep2}
    q_{\mu}(t) = \lambda {e^{-\lambda^2/4t} \o \sqrt{\pi t} }
    \int_0^\infty (e^{-\kappa/4t} - 1 + \kappa/4t)w_\lambda(v)dv.
    \end{equation}
\end{thm}
\begin{proof}
 The proof relies on application of the formula \pref{expintegr1} for the function
 \begin{equation*}
 \hat{\varphi}(r^2)= x^{\mu} { K_{\mu}(xr) \o K_{\mu}(r)}\/.
 \end{equation*}
 Note that by \pref{Flambda0} and \pref{laplace01}  it follows that
 \begin{eqnarray*}
    rx^{\mu} { K_{\mu}(xr) \o K_{\mu}(r)} &=& e^{-\lambda r} \lambda F_{\lambda}(r)
    +  e^{-\lambda r}  x^{\mu-1/2} Q(r) \\
    &=& e^{-\lambda r}  x^{\mu-1/2} r +  e^{-\lambda r} \lambda F_{\lambda}(r) -
     e^{-\lambda r}x^{\mu-1/2} (\mu^2 - 1/4) \lambda/2x \\
    &=&  e^{-\lambda r}  x^{\mu-1/2} r +
    \lambda e^{-\lambda r} \int_0^\infty \(e^{-rv} - 1\) w_\lambda(v) dv.
  \end{eqnarray*}
 We recall that $Q(r)= r - (\mu^2-1/4) \lambda/2x $ and
 $F_{\lambda}(r)= \int_0^\infty e^{-rv} w_{\lambda}(v) dv$.
 To simplify the proof we introduce new notation. Namely, denote
  \begin{equation*}
   w_{\lambda}^{\#}(v)=-\int_v^\infty  w_{\lambda}(v) dv,
   \qquad  \text{that is} \qquad    {d w_{\lambda}^{\#}(v) \o dv}= w_{\lambda}(v) .
 \end{equation*}
 Then, by integration by parts we obtain
 \begin{eqnarray*}
  \int_0^\infty \(e^{-rv}-1\)  w_{\lambda}(v) dv &=& \(e^{-rv}-1\)
   w_{\lambda}^{\#}(v) |_{0}^{\infty}
  + r \int_0^\infty e^{-rv} w_{\lambda}^{\#}(v) dv \\
  &=&  r \int_0^\infty e^{-rv} w_{\lambda}^{\#}(v) dv .
    \end{eqnarray*}
 Thus, the left-hand side of the formula \pref{expintegr1}
 is now of the form
 \begin{eqnarray*}
 &{}& \int_0^{\infty}  rx^{\mu} { K_{\mu}(xr) \o K_{\mu}(r)} \cos(ry)/r dr \\
 &=&  x^{\mu-1/2} \int_0^{\infty} e^{-\lambda r}  \cos(ry) dr \\
 &+&   \int_0^{\infty} \lambda e^{-\lambda r} \{\int_0^\infty
    \( e^{-rv} -1\) w_{\lambda}(v) dv \} \cos(ry)/r dr\\
 &=&  x^{\mu-1/2} \int_0^{\infty} e^{-\lambda r}  \cos(ry) dr \\
 &+&  \int_0^{\infty} \lambda e^{-\lambda r} \{\int_0^\infty
     e^{-rv}  w_{\lambda}^{\#}(v) dv \} \cos(ry) dr\\
 &=& H_{\mu}^{(1)}(y) +  H_{\mu}^{(2)}(y) =  H_{\mu}(y).
  \end{eqnarray*}
   Using the standard formula for the Laplace transform of the cosine function
 we obtain
  \begin{equation*}
  H_{\mu}^{(1)}(y)
   = x^{\mu-1/2}
       {\lambda   \o \lambda^2+y^2}\/,
  \end{equation*}
   and
  \begin{equation*}
   H_{\mu}^{(2)}(y)
   = \lambda
   \int_0^\infty {(\lambda +v)
    w_\lambda^{\#}(v) \o (\lambda+v)^2+y^2} dv.
    \end{equation*}
 By the form of $ H_{\mu}^{(1)}(y) $, \pref{kernel} (applied for $s=1$),
  \pref{laplace01} and \pref{laplace02}
 we obtain
 \begin{equation*}
  H_{\mu}^{(1)}(y) ={ x^{\mu-1/2} \lambda \o 2^{2}}
 \int_0^\infty e^{-y^2/4t} e^{-\lambda^2/4t} {dt \o t^{2}}.
  \end{equation*}
  Analogously, for $H_{\mu}^{(2)}$ we obtain
  \begin{eqnarray*}
  H_{\mu}^{(2)}(y) &=&
  { \lambda \o 2^{2}}
      \int_0^\infty e^{-y^2/4t} \(\int_0^\infty (\lambda +v) w_\lambda^{\#}(v) e^{-(\lambda+v)^2/4t} dv \)
        {dt \o t^{2}} \\
   &=&
    { \lambda \o 2}
        \int_0^\infty e^{-y^2/4t} \(\int_0^\infty {d(- e^{-(\lambda+v)^2/4t}) \o dv}
         w_\lambda^{\#}(v) dv \)
          {dt \o t}. \\
  \end{eqnarray*}
 Observe now that
 \begin{eqnarray*}
  \int_0^\infty {d(- e^{-(\lambda+v)^2/4t}) \o dv} w_\lambda^{\#}(v) dv &=&
 - e^{-(\lambda+v)^2/4t} w_\lambda^{\#}(v) |_0^\infty +
 \int_0^\infty  e^{-(\lambda+v)^2/4t} {d w_\lambda^{\#}(v) \o dv} dv \\
 &=& e^{-\lambda^2/4t}  w_\lambda^{\#}(0) +
  \int_0^\infty  e^{-(\lambda+v)^2/4t} w_\lambda(v) dv.
  \end{eqnarray*}
  We also have
  \begin{equation*}
  w_\lambda^{\#}(0) = -\int_0^\infty {d w_\lambda^{\#}(v) \o dv} dv = -\int_0^\infty
  w_\lambda(v) dv.
  \end{equation*}
 By the above identities and the form of $H_{\mu}^{(2)}$ we obtain

 \begin{equation*}
   H_{\mu}^{(2)}(y) =
  { \lambda \o 2}
      \int_0^\infty e^{-y^2/4t} \(\int_0^\infty ( e^{-\kappa/4t} -1) w_\lambda(v)  dv \)
        {dt \o t}.
  \end{equation*}

 Combining the above identities we obtain that the left-hand side of \pref{expintegr1}
 takes the form
 \begin{equation*}
   H_{\mu}(y) = {\lambda \o 2} \int_0^\infty  e^{-y^2/4t}  e^{-\lambda^2/4t}
   \(x^{\mu-1/2}/2t + \int_0^\infty  ( e^{-\kappa/4t} -1)  w_\lambda(v) dv \) {dt \o t}.
   \end{equation*}
 Taking into account the right-hand side of \pref{expintegr1} and continuity of $q_{\mu}$
 with respect to $t$ (see properties of $w_\lambda$ below), we obtain \pref{rep1}.
 When $\mu>1/2$, then using \pref{laplace02} we obtain \pref{rep2}.

 \end{proof}

Below we give a description of the function $w_\lambda$. We rely here on
results contained in \cite{BGS}.
The formulas depend on the zeros of the function $K_\mu(z)$.

Even if in general the values of these zeros are not
given explicitly,  we are able to prove some important properties
(as boundedeness or asymptotics) of $w_\lambda$, which are essential
 in applications. Moreover, for some values of $\mu$ we provide
 explicit formulas as well (see Corollary \ref{examples}).

The function $K_\mu(z)$ extends  to an entire function
when $\mu-1/2$ is an integer and has a holomorphic extension to
$\C\setminus (-\infty,0]$ when $\mu-1/2$ is not an integer.
Denote the set of
zeros of the function $K_\mu(z)$ by
$Z=\{z_1,...,z_{k_\mu}\} $. We give some information about
these zeros (cf. \cite{E}, p. 62) needed in the sequel.
Recall that $k_\mu=\mu-1/2$
when $\mu-1/2 \in \N$. For $\mu-1/2 \notin \N$,
$k_\mu$ is the even number
closest to $\mu-1/2$. In particular, for $0\leq \mu<3/2$ we have $k_\mu
=0$; for $\mu=2$ and $3$ we have $k_\mu=2$. The functions
$K_\mu$ and $K_{\mu-1}$ have no common zeros.

 As in \cite{BGS}, we need an additional notation to describe the function $w_\lambda$.
 Define for $\mu>0$
\begin{equation} \label{w_1}
w_{1,\/\lambda}(v)= -{x^\mu \o \lambda} \sum_{i=1}^{k_\mu}
 {z_i e^{\lambda z_i} K_\mu(xz_i) \o K_{\mu-1}(z_i)} \/ e^{z_i v}\/.
\end{equation}

When $\mu+1/2 \notin \N$  and $\mu\ge0$ we define
\begin{equation} \label{w_2}
 w_{2,\/\lambda}(v) =-\cos(\pi\mu) {x^\mu \o \lambda}
 \int_0^\infty {
   I_\mu\(xu\)K_\mu(u) - I_\mu(u)K_\mu\(xu\)\o
\cos^2(\pi\mu) K_\mu^2(u)+(\pi I_\mu(u)+\sin(\pi \mu) K_\mu(u))^2} \/
 e^{-\lambda u} e^{-vu} \/u du\/.
\end{equation}

We now formulate our representation theorem for the function
$w_\lambda$. The proof of the main part is the same as in \cite{BGS} and is omitted;
we only show asymptotic properties of the function $w_{2,\lambda}$. For $\mu=0$ this is
new; behaviour for $\mu=(n-1)/2$ was shown in \cite{BGS}. Nevertheless, we present here
a new and unified proof based on tauberian theorems.

\begin{thm}\label{formula}
 In the case $\mu-1/2 \in \N$
\begin{equation*}
w_\lambda(v) = w_{1,\/ \lambda}(v);
\end{equation*}
while, in the case when  $\mu-1/2 \notin \N$
\begin{equation*} 
w_\lambda(v) = w_{1,\/ \lambda}(v) +  w_{2,\/ \lambda}(v)\/.
\end{equation*}
Moreover, we have $\sup_{v\geq 0}|w_{\lambda}(v)| <\infty\/, $ and
$$
- \cos(\pi \mu)w_{2,\lambda}(v)\geq 0, \ \ \ v\geq 0,\ \ \ \  (\mu-1/2 \notin \N);
$$
\begin{equation*}
   \int_0^\infty v^k |w_{1,\lambda}(v)| dv <\infty, ~~ k=1, 2,\ldots \/;
\end{equation*}
\begin{equation*}
 \lim_{v \to \infty} v^k w_{1\/,\/\lambda}(v) = 0,\quad  k=1,2,\ldots \/;
\end{equation*}
\begin{equation*}
 \lim_{v \to \infty} v^{2\mu+2} w_{2\/,\/\lambda}(v) = {-\cos(\pi \mu) \Gamma(2\mu+2)
   \o 2^{2\mu-3} \Gamma(\mu) \Gamma(\mu +1)}
   {x^{2\mu}-1 \o \lambda},\/\ \ \   (\mu-1/2 \notin \N,\ \mu>0)\/;
\end{equation*}
 \begin{equation*}
\lim_{v \to \infty}(v \log v)^2 w_{2,\/\lambda}(v)= -{1\o 2\lambda}\log x ,\quad \quad
 \text{for $\mu=0$}\/.
 \end{equation*}
\end{thm}
\begin{proof} Denote
\begin{equation} \label{h_2}
 h_{\mu,\/\lambda}(u) =
   {I_\mu\(xu\)K_\mu(u) - I_\mu(u)K_\mu\(xu\)\o
\cos^2(\pi\mu) K_\mu^2(u)+(\pi I_\mu(u)+\sin(\pi \mu) K_\mu(u))^2}  e^{-\lambda u}u\/.
\end{equation}
Observe that the function $ h_{\mu,\/\lambda}(u)$ is non-negative
($I_\mu(u)/K_\mu(u)$ is increasing for $u>0$)
and $w_{2,\/\lambda}(v)$ is the Laplace transform of
 $-\cos(\pi\mu) {x^\mu \o \lambda}h_{\mu,\/\lambda}$ at $v$.
We  claim that for $u\to 0+$ we have the following asymptotics for $h_{\mu,\/\lambda}$:
\begin{equation}\label{asym}
h_{\mu,\/\lambda}(u)\approx\begin{cases}  x^\mu\frac {c_\mu}{c'_\mu}(1-x^{-2\mu})
 u^{2\mu+1}\/,&\text{for $\mu>0$}\/,\\
    {u \/ \log x  (\log u)^{-2}}, &\text{for $\mu=0$}\/.
\end{cases}
\end{equation}
Applying Karamata's Tauberian theorem (see, e.g. \cite{Fe}) we obtain
that the  asymptotic behaviour of $w_{2,\/\lambda}(v)$ is
$$\lim_{v \to \infty}w_{2,\/\lambda}(v)\left(\int_0^{1/v}
 h_{\mu,\/\lambda}(u)\right)^{-1}= {-\cos(\pi \mu)\o \lambda}\Gamma(2+2\mu)x^\mu.$$
This together with  \pref{asym} implies:
\begin{eqnarray*}
&{}&
 \lim_{v \to \infty} v^{2\mu+2} w_{2\/,\/\lambda}(v) = {-\cos(\pi \mu) \Gamma(2\mu+2)
   \o 2^{2\mu-3} \Gamma(\mu) \Gamma(\mu +1)}
   {x^{2\mu}-1 \o \lambda},\/\ \ \   (\mu-1/2 \notin \N,\ \mu>0)\/, \\
&{}&
\lim_{v \to \infty}(v \log v)^2 w_{2,\/\lambda}(v)= -{1\o 2\lambda}\log x ,
\quad \quad
  \mu=0 \/.
\end{eqnarray*}
To prove \pref{asym} we apply the following asymptotics.
When $u \to 0$ we have:
\begin{eqnarray}
I_\mu(u) &\approx& c_\mu u^\mu, \quad
K_\mu(u) \approx c'_\mu u^{-\mu}\/;\\ \label{kizero}
I_0(u) &=& 1+o(1), \quad \quad
K_0(u) = \log (2/u)I_0(u)+\psi(1)+o(1)\/, \label{kizero_1}
\end{eqnarray}
with $c_\mu=2^{-\mu}/\Gamma(\mu+1)$,
$c'_\mu=2^{\mu-1}\Gamma(\mu)$ and where $\psi$ is the Euler function. Then
\begin{eqnarray*} \label{h_2_1}
 e^{\lambda u}u^{-1}h_{\mu,\/\lambda}(u) &=&
   {I_\mu\(xu\)K_\mu(u)\o
 K_\mu^2(u)+(\pi I_\mu(u))^2+2\pi\sin(\pi \mu) K_\mu(u)I_\mu(u)}
 \left( 1-{I_\mu(u)K_\mu\(xu\)\o I_\mu\(xu\)K_\mu(u)}\right)\\
 &\approx& {I_\mu\(xu\)\o  K_\mu(u)}
  \left( 1-{I_\mu(u)K_\mu\(xu\)\o I_\mu\(xu\)K_\mu(u)}\right)\\
 &\approx&
 \begin{cases}
    u^{2\mu}x^\mu\frac {c_\mu}{c'_\mu}(1-x^{-2\mu}), &\text{for $\mu>0$}\/,\\
    {\log x \o (\log u)^2}, &\text{for $\mu=0$}\/.
\end{cases}
\end{eqnarray*}
\end{proof}

{\bf Examples.}

  To illustrate representation theory developed so far  we write down explicit
   integral formulas  for the density
 $q_\mu$ in some special cases of $\mu$. All formulas appearing here follow directly
  from Theorems
 \ref{rep}  and \ref{formula}. If $0 \leq \mu<3/2$  then
 $w_\lambda=w_{2,\lambda}$ and functions  $w_\lambda$ have constant sign. For $\mu=0$
 and $\mu=1$ the function $w_\lambda$ has simpler form, which we exhibit here.
 If $\mu=1/2$ then the function $F_\lambda=0$ and $q_{\mu}(t)$
 reduces to the standard $1/2$-stable subordinator. For $\mu+1/2 \in \N$
 we have, in turn that $w_\lambda=w_{1,\lambda}$ and the form of $w_{1,\lambda}$
 can be computed calculating residues of simple rational functions (see the formula
 for $w_{1,\lambda}$ or calculations in \cite{BGS}). Again, we write the explicit form
 of $w_\lambda$  for $\mu=3/2$ and  $\mu=5/2$.

\begin{cor}\label{examples}
If $\mu=0$ then
  \begin{equation*}
   - w_\lambda(v)={1 \o \lambda } \int_0^\infty { I_0(xu)
      K_0(u) - K_0(xu) I_0(u) \o K_0^2(u)+\pi^2 I_0^2(u)} e^{-u
      \lambda} e^{-vu}u\/du\/,
  \end{equation*}
 and
  \begin{equation*}
    q_{\mu}(t) = \lambda {e^{-\lambda^2/4t} \o \sqrt{\pi t} }
   \((\lambda +1)^{-1/2}/2t +
    \int_0^\infty \(1-e^{-\kappa/4t} \)(- w_\lambda(v))\/dv \)\/.
\end{equation*}

 If $\mu=1/2$ then
\begin{equation*}
q_{\mu}(t) = \lambda {e^{-\lambda^2/4t} \o 2 \sqrt{\pi t^3} }\/.
\end{equation*}

If $\mu=1$ then
  \begin{equation*}
    w_\lambda(v)={\lambda+1 \o \lambda } \int_0^\infty { I_1(xu)
      K_1(u) - K_1(xu) I_1(u) \o K_1^2(u)+\pi^2 I_1^2(u)} e^{-u
      \lambda} e^{-vu}u\/du\/,
  \end{equation*}
 and
 \begin{equation*}
     q_{\mu}(t) = \lambda {e^{-\lambda^2/4t} \o \sqrt{\pi t} }
     \int_0^\infty (e^{-\kappa/4t} - 1 + \kappa/4t)w_\lambda(v)dv.
     \end{equation*}

   If $\mu=3/2$ then $w_\lambda(v)= e^{-v}$ and
  \begin{equation*}
   q_{\mu}(t) = \lambda {e^{-\lambda^2/4t} \o \sqrt{\pi t} }
      \int_0^\infty (e^{-\kappa/4t} - 1 + \kappa/4t)e^{-v}  dv.
     \end{equation*}

  If $\mu=5/2$ then
  \begin{equation*}
    w_\lambda(v)= 3 e^{-3v/2}
[(2\lambda
    +1)\cos(\sqrt{3} v/2) +\sqrt{3}  \sin(\sqrt{3} v/2)] \/,
  \end{equation*}
 and
 \begin{equation*}
      q_{\mu}(t) = \lambda {e^{-\lambda^2/4t} \o \sqrt{\pi t} }
      \int_0^\infty (e^{-\kappa/4t} - 1 + \kappa/4t)w_\lambda(v)dv.
      \end{equation*}

 \end{cor}

\section{Asymptotic behaviour of $A(\tau)$}\label{asymp}

In this section we prove the following
  \begin{thm} \label{Ainfty}
   The density $q_\mu$ of $A(\tau)$ satisfies:
  $$ \lim_{t\to\infty} t^{\mu +1}q_\mu(t)=C_{\mu},\quad \text{if $\mu>0$}\/,$$
  $$ \lim_{t\to\infty} (\log t)^2tq_\mu(t)=C_0,\quad \text{if $\mu=0$}\/,$$
  for some positive $C_\mu$\/.
   \end{thm}
\begin{lem}\label{estimate}
Let $\mu>0$. There exists a constant $C >1$ such that
$$C^{-1} t^{-\mu}
\le P(A(\tau)>t)\le C t^{-\mu},\quad t>1.$$
\end{lem}
\begin{proof}
 From (\ref{basic}) we infer  that $(x>1, a=1)$:
$$ P(A(\infty)>t)\le P(A(\tau)>(1-x^{-1})tx^2)+P(A(\infty)>x t), $$
 which implies
 \begin{equation}\label{dom1} P(t\le A(\infty)\le xt )\le P(A(\tau)>(x-1)xt).
 \end{equation}
 Moreover,
 \begin{equation} \label{dom2} P(A(\tau)>x^2t)\le P(A(\infty)>t).\end{equation}
Now the lemma follows from (\ref{dom1}) and  (\ref{dom2}) since
 $P(A(\infty)>t)\approx c t^{-\mu}$ for some $c>0$ by (\ref{subord}).
 \end{proof}

 Recall that $\kappa=(\lambda+v)^2-\lambda^2$.

  \begin{lem}\label{katymoment0}
Let $m \in \N$ be such that $2\le m\le\mu +1/2$.
 Then
\begin{equation} \label{cm}
\lim_{t\to\infty} t^{m} \int_0^\infty w_\lambda(v)
 (e^{-\kappa/4t}-\sum_{0\leq j\leq m-1} (-1)^j{1\o j!}
 \left(\frac\kappa {4t}\right)^j)dv =
  \frac{(-1)^m}{4^m m!}\int_0^\infty\kappa^m w_\lambda(v)\ dv=C_m\/.
 \end{equation}
 Moreover,  $C_m=0$, for $2\le m< \mu+1/2$, and, in the case $m=\mu+1/2\in \N$,
 we have $C_m>0$.

\end{lem}
\begin{proof} Denote
$$\psi(t,m)=t^{m}  w_\lambda(v)
 (e^{-\kappa/4t}-\sum_{0\leq j\leq m-1} (-1)^j{1\o j!}
  \left(\frac\kappa {4t}\right)^j).$$
 By elementary calculations
 $$
 \left|e^{-\kappa/4t}-\sum_{0\leq j\leq m-1} (-1)^j{1\o j!}
 \left(\frac\kappa {4t}\right)^j \right|\le {1 \o m!} \(\frac\kappa {4t}\)^{m}.
 $$
 Hence
 $$|\psi(t,m)|\le |w_\lambda(v)|{\kappa^{m}\o m!}.$$
 Under the assumption on $m$ the function $|w_\lambda(v)|\kappa^{m}$ is integrable so
 the formula \pref{cm} follows from the  bounded convergence theorem.

Suppose that there exists $C_m \neq 0$, with $m$ having properties as above, and
 denote
 \begin{equation*}
 m_0 =\inf \{ m\in\N: 2\le m \leq \mu+1/2, C_m \neq 0\} \/.
 \end{equation*}
  Then from the first part of the proof we have
 $$ \lim_{t\to\infty} t^{m_0+1/2}q_\mu(t)=C_{m_0}>0,$$\textbf{}
 which implies that
 $$\lim_{t \to\infty} t^{m_0-1/2}P(A(\tau)>t)=C_{m_0}(m_0-1/2)^{-1}. $$
   From Lemma \ref{estimate} we infer that $C_{m_0}>0$ if and only if
  $m_0-1/2=\mu$. In particular, we then have $\mu+1/2\in \N$.
  We also obtained that if $m < \mu+1/2$ then $C_m=0$.
 This completes the proof of the lemma.
 \end{proof}

  {\bf Remark.}
  The above lemma yields, in particular, that for $m\in\N,\quad 2\leq m < \mu+1/2$
    \begin{equation*}
    \int_0^\infty \kappa^m w_\lambda(v)\/ dv =0\/.
    \end{equation*}
    Thus, the representation formula for the density $q_\mu$
    can be written for $\mu \geq 1/2$ as follows:
  \begin{equation}\label{rep2n}
      q_{\mu}(t) = \lambda {e^{-\lambda^2/4t} \o \sqrt{\pi t} }
      \int_0^\infty (e^{-\kappa/4t} - \sum_{0\leq j\leq l} (-1)^j{1\o j!}
      \left(\frac\kappa {4t}\right)^j )w_\lambda(v)dv,
      \end{equation}
  where $l=[\mu+1/2]$, if $\mu-1/2\notin \N$, and $l=\mu-1/2$ otherwise.

 We now prove our theorem.

 {\bf Proof of Theorem \ref{Ainfty}}

 \begin{proof}
  For $\mu=1/2$ the density
 $q_\mu$ has a particularly simple form (see Corollary \ref{examples}) and
 the theorem clearly holds true.
 Hence, we assume throughout the remainder of the proof that $\mu \neq 1/2$.
 Next, if $\mu-1/2 \in \N$ then the Remark above together with Lemma
 \ref{katymoment0} yield our theorem at once.

 Thus, we assume for what follows that $\mu-1/2 \notin \N$ and let
 $l= [\mu+1/2]$.
 Denote
 $$I(t)= \int_0^\infty w_\lambda(v) t^{\mu +1/2}
 (e^{-\kappa/4t}-\sum_{0\leq j\leq l} (-1)^j{1\o j!}
  \left(\frac\kappa {4t}\right)^j)dv.$$
  We prove that
  \begin{equation} \label{asIt}
  \lim_{t \to \infty} I(t) =C>0\/, \quad \text{if} \quad \mu>0\/;
  \end{equation}
 and
  \begin{equation} \label{aslogIt}
  \lim_{t \to \infty}(\log t)^2 I(t)=C>0\/, \quad \text{if} \quad \mu=0\/.
  \end{equation}
Applying change of variable $\kappa=4st$ we obtain
  $$
  v=\sqrt{4st+\lambda^2}-\lambda, \quad dv= \frac {2t}{\sqrt{4st+\lambda^2}}ds,
   $$
   so
   $$I(t)= \int_0^\infty \psi_\lambda(s,t)
 (e^{-s}-\sum_{0\leq j\leq l} (-1)^j{s^j\o j!})
  ds,$$
where
\begin{equation*}
\psi_\lambda(s,t)= w_\lambda(\sqrt{4st+\lambda^2}-\lambda)
    \frac {2t^{\mu +3/2}}{\sqrt{4st+\lambda^2}}\/.
\end{equation*}

  We claim that for $\mu>0$ and $t,s\in R^+$ there is a constant $C>0$ such that
 \begin{equation} \label{psi1}
 |\psi_\lambda(s,t)| \leq C s^{-\mu -3/2}\/.
 \end{equation}
 For $\mu=0$ and  $t,s\in R^+$  our claim is:
 \begin{equation} \label{psi2}
 (\log t)^2 |\psi_\lambda(s,t)| \leq C\max\{1,(\log s)^2\} s^{-3/2}\/,
  \end{equation}
 where $C>0$.
  The above claims prove the relations \pref{asIt} and \pref{aslogIt}. Indeed,
 consider first the case $\mu>0$. Then the absolute value of the expression under the
 integral $I(t)$ can be
   estimated by the integrable function
   $$C \left|e^{-s}-\sum_{0\leq j\leq l} (-1)^j{s^j\o j!}\right|s^{-\mu-3/2}\le C \min\{1,s\}s^l s^{-\mu-3/2}=
   C \min\{1,s\} s^{-\delta-1},$$
   where $\delta=\mu+1/2-[\mu+1/2]$, $0<\delta<1$. Then the proof of \pref{asIt}
   is concluded by passing
   $t \to \infty$ and using the asymptotics of $w_\lambda$ (see Theorem \ref{formula}).

Now, consider the case $\mu=0$. Observe that the absolute value of the
integrand in $(\log t)^2 I(t)$ is estimated by the integrable function
 $$ C\max\{1,(\log s)^2\} s^{-3/2}|e^{-s}-1|\/,
so
$$
 $$\lim_ {t\to \infty} (\log t)^2 I(t)={\log x\o 2\lambda}\int_0^\infty \( 1-e^{-s}\)s^{-3/2}ds.$$
Here we take into account
 $$
  (\log t)^2 \psi_\lambda(s,t)
  \approx  -
  {\log x\o 2\lambda}{(\log t)^2\o s^{3/2}(\log st)^2}\approx -
  {\log x\o 2\lambda s^{3/2}}, \quad t \to \infty,$$
   by  using the asymptotics of $w_\lambda$ (see Theorem \ref{formula}).


 Now we can conclude the proof of the asymptotic behaviour of $q_\mu(t)$ in the case
  $\mu -1/2\notin \N$. Note that for $0<\mu<1/2$ we have
 \begin{eqnarray*}
    t^{\mu+1} q_{\mu}(t) &=&\lambda {e^{-\lambda^2/4t} \o \sqrt{\pi t} }
    \(t^{\mu} x^{\mu-1/2}/2 +
     t^{\mu+1}\int_0^\infty \(e^{-\kappa/4t} - 1\) w_\lambda(v)dv \)\\
     &=&\lambda {e^{-\lambda^2/4t} \o \sqrt{\pi} }
    \(t^{\mu -1/2} x^{\mu-1/2}/2 + I(t)\)\rightarrow \lambda C/\sqrt{\pi},
     \end{eqnarray*}
 where $C$ is the  constant from the right-hand side of the formula \pref{asIt}.
 The same argument for $\mu=0$ shows that
 $$(\log t)^2 t q_{\mu}(t)\rightarrow \lambda C/\sqrt{\pi} \/,$$
 where $C$ is the  constant from the right-hand side of the formula \pref{aslogIt}.
 For $\mu>1/2$ the asymptotics of $q_{\mu}(t) $ directly follows from
 the formula \pref{asIt}.
  \end{proof}
 We now justify our claims \pref{psi1} and \pref{psi2}. We use
 notation as introduced in the proof of the theorem.
 \begin{lem} \label{claim}
 If  $\mu>0$ there is $c>0$ such that
  \begin{equation*}
  |\psi_\lambda(s,t)| \leq c s^{-\mu -3/2}\/.
  \end{equation*}

   For $\mu=0$ there is $C>0$ such that
  \begin{equation*}
  (\log t)^2 |\psi_\lambda(s,t)| \leq C \max\{1,(\log s)^2\} s^{ -3/2}\/.
   \end{equation*}
 \end{lem}
 \begin{proof}
  We begin with \pref{psi1} first.
 Let $\mu>0$. Using boundedness of the function $w_\lambda$ we can estimate for
  $ts\le 1$:
    $$
    |\psi_\lambda(s,t)|
   \le \sup_{v\ge 0}|w_\lambda(v)|
     \frac {t^{\mu +3/2}}{\sqrt{4st+\lambda^2}}
    \le C s^{-\mu-3/2}.
    $$
   For $ts\ge 1$ we use the asymptotics of $w_\lambda$ at $\infty$
   (see Theorem \ref{formula})
    to arrive at
     $$
     |\psi_\lambda(s,t)|
    \le C {t^{\mu +3/2}}(st)^{-\mu-3/2}= C s^{-\mu-3/2}\/.
    $$
   For $\mu=0$ we estimate our expression as follows.
 Again using boundedness of the function $w_\lambda$
      we estimate for $ts\le 2$:
     $$
     |\psi_\lambda(s,t)|
     \le \sup_{v\ge 0}|w_\lambda(v)| \frac { t^{3/2}}{\sqrt{4st+\lambda^2}}
     \le C  s^{-3/2}.
     $$
       For $2\le ts$ we use the asymptotics of $w_\lambda$ to get
      $$
   (\log t)^2 |\psi_\lambda(s,t)|
    \le C {t^{3/2}}(st)^{-3/2}{(\log t)^2\o (\log st)^2}
    \le \begin{cases}
    C s^{-3/2}(\log t)^2, & 2\le ts\le \sqrt{t}\/,\\
    C s^{-3/2},&   ts\ge \sqrt{t}\/.
    \end{cases}
     $$
Next observe that for $2\le ts\le \sqrt{t}$ we have $(\log t)^2\le 4(\log s)^2$.
 \end{proof}

\section{Hiperbolic Brownian motion with drift}

Consider the half-space model of the $n$-dimensional real hyperbolic space
\begin{equation*}
\H^n = \{ (x_1,\ldots,x_{n-1},x_n)\in \R^{n-1}\times \R:\; x_n>0 \}  .
\end{equation*}
The Riemannian metric, the volume element and the Laplace-Beltrami
operator are given by
\begin{equation*}
  ds^2 = {dx_1^2+...+dx_{n-1}^2+dx_n^2 \o x_n^2},
\end{equation*}
\begin{equation*}
  dV={dx_1...dx_{n-1}dx_n \o x_n^n},
\end{equation*}
\begin{equation*}
  \Delta= x_n^2 (\sum_{i=1}^{n} \partial_i^2)- (n-2)x_n\partial_n,
\end{equation*}
respectively (here $\partial_i={\partial \o\partial x_i}$, $i=1,...,n$).
For $\mu \ge 0$ let $\alpha = 2\mu-n+1$. We also introduce the operator:
\begin{equation*}
  \Delta_\mu =\Delta - \alpha x_n\partial_n=  x_n^2 (\sum_{i=1}^{n} \partial_i^2)- (2\mu-1)x_n\partial_n.
\end{equation*}
Let $(B_i(t))_{i=1...n}$ be a family of independent classical
Brownian motions on $\R$ with the generator ${d^2\o dx^2}$ (and not
${1\o 2}{d^2 \o dx^2}$) i.e.  the variance $E^0B_i^2(t) = 2t$.
Then the Brownian motion on $\H^n$, $X=(X_i)_{i=1...n}$ 
can be described by the following system of stochastic differential
equations
\begin{equation*}
  \left\{
    \begin{array}{ccc}
      dX_1(t) & = & X_n(t)dB_1(t) \\
      dX_2(t) & = & X_n(t)dB_2(t) \\
      . & . &.  \\
      dX_n(t) & = & X_n(t)dB_n(t) - (n-2)X_n(t)dt.
    \end{array} \right.
\end{equation*}
More generally, if we replace $n-2$ by $\alpha + n-2=2\mu-1$ that the corresponding process will be called the Brownian motion on $\H^n$ with drift $\alpha$.
By the It\^o formula one verifies that the generator of the solution
of this system is $\Delta_\mu$. Moreover, it can be easily checked that
the solution is given by
\begin{equation*}
  \left\{
    \begin{array}{rcl}
      X_1(t) & = & X_1(0)+\int_0^t X_n(t) dB_1(s) \\
      X_2(t) & = & X_2(0)+\int_0^t X_n(t) dB_2(s)\\
      . & . &.  \\
      X_n(t) & = &X_n(0)\exp( B_n(t)-2\mu t).
    \end{array} \right.
\end{equation*}

Define the projection $\tilde {} : \R^n\ni u=(u_1,...,u_n) \to \tilde
u=(u_1,...,u_{n-1})\in \R^{n-1}$.  In particular,  $\tilde X(t) =
(X_1(t),...,X_{n-1}(t))$. From the representation above one may easily verify
( e.g. by computing  Fourier transforms) that
\begin{equation}\label{repAt}
 \tilde X(t)\stackrel{d}= \tilde X(0) +\tilde B(\int_0^t X^2_n(s) ds), \quad t\ge 0.
 \end{equation}
where $\tilde B(t)= (B_1(t),...,B_{n-1}(t))$ is $n-1$ dimensional  Brownian
 motion independent of the process $X_n(t)$.
Consider a half-space $D=\{ u\in \H^n:\; u_n>a \}$ for some fixed
$a>0$. To simplify the notation we choose $a=1$. Define
\begin{equation*}
  \tau=\inf\{ t\geq 0:\; X(t)\notin D\} =
\inf\{ t\geq 0:\; X_n(t) = 1\}.
\end{equation*}
By $P_1(u,y)$, $u=(u_1,u_2,...,u_n)\in D$,
$y=(y_1,y_2,...,y_{n-1},1)\in \partial D$ we denote the
Poisson kernel of $D$, ie. the density of the distribution  of $X(\tau)$
starting at $u$ (since $X_n(\tau)=1$ it is enough to consider
the distribution of $\tilde X(\tau)$).
From \pref{repAt} it is obvious that
$$\tilde X(\tau)\stackrel{d}=\tilde u + \tilde B(A(\tau)),$$
where
the functional    $A(\tau)$ (starting from $u_n>1$) is independent of $\tilde B(t)$.
For further considerations we may take  $\tilde u=0$ and $u_n=x>1$,
so the starting point of the $n$ dimensional process  $X(\cdot)$
is $(0,\dots,0, x)\in D$. Since $A(\tau)$ and $\tilde B(t)$ are independent,
we have the following representation of the Poisson kernel:
\begin{cor}\label{Poisson}
$$P_1(x,y)= {1 \o (4\pi)^{(n-1)/2}} \int_0^\infty
e^{-|y|^2/4t} q_{\mu}(t){dt \o t^{(n-1)/2}}  . $$
\end{cor}
Observe that for $\mu=1/2$ the functional $A(\tau)$ has the standard asymmetric
 $1/2$-stable  distribution and the Poisson kernel is exactly $n-1$
 dimensional Cauchy density, so in what follows, we consider $\mu\neq 1/2$.

For the sake of simplicity we write $|y|$ as $\rho$.
Taking into account the formula \pref{subord} and \pref{Pkernel} we obtain
the following representation of the Poisson kernel of the set $D$
\begin{thm} \label{Pkernrepr}
For $0\leq \mu < 1/2$ we have
\begin{equation*}
P_1(x,y)= {\Gamma({n \o 2}-1) \o 2 \pi^{n/2}} {\lambda \o (\lambda^2+\rho^2)^{n/2}}
\[(n-2)(\lambda+1)^{\mu -1/2} -
 \int_0^\infty { w_\lambda(v) L^{\#}(\lambda,\rho,v) \/dv \o
((\lambda +v)^2 +\rho^2)^{{n \o 2}-1}} \]\/,
\end{equation*}
where $L^{\#}(\lambda,\rho,v)$ is the following function
\begin{equation*}
 L^{\#}(\lambda,\rho,v) =
 (\lambda^2+\rho^2) [((\lambda+v)^2+\rho^2)^{n/2-1}-(\lambda^2+\rho^2)^{n/2-1}].
\end{equation*}

For $\mu>1/2$ we obtain
\begin{equation*}
P_1(x,y)= {\Gamma({n \o 2}-1) \o 2 \pi^{n/2}} {\lambda \o (\lambda^2+\rho^2)^{n/2}}
\int_0^\infty { w_\lambda(v) L(\lambda,\rho,v) \/dv \o
((\lambda +v)^2 +\rho^2)^{{n \o 2}-1}}\/,
\end{equation*}
with $L(\lambda,\rho,v)$ defined by
\begin{eqnarray*}
 &{}&L(\lambda,\rho,v) \\
 &=& (n/2-1)((\lambda+v)^2-\lambda^2)((\lambda+v)^2+\rho^2)^{n/2-1}
 - (\lambda^2+\rho^2) [((\lambda+v)^2+\rho^2)^{n/2-1}-(\lambda^2+\rho^2)^{n/2-1}].
\end{eqnarray*}
\end{thm}
\begin{proof}
By the formula \pref{kernel} we obtain
\begin{equation*}
\int_0^\infty e^{-\rho^2/4t} e^{-\lambda^2/4t} {dt \o t^{1+s}}
 = { 2^{2s}\Gamma(s) \o (\lambda^2+ \rho^2)^{s}}.
\end{equation*}
Applying the above formula (with different constants) for three terms appearing
in the representation of $q_\mu$
we obtain,

\begin{equation*}
\int_0^\infty e^{-\rho^2/4t} e^{-\lambda^2/4t} {dt \o t^{1+n/2}}
 = { 2^n \Gamma({n \o 2}) \o (\lambda^2+\rho^2)^{n/2}}\/,
\end{equation*}
 and
\begin{eqnarray*}
&{}& \int_0^\infty e^{-\rho^2/4t} e^{-\lambda^2/4t}
\int_0^\infty e^{-\kappa/4t} w_\lambda(v) {dt \o t^{n/2}} \\
&=& \Gamma({n \o 2}-1) 2^{n-2} \int_0^\infty
{ w_\lambda(v) dv \o ((\lambda+v)^2+\rho^2)^{{n \o 2}-1}}\/,
\end{eqnarray*}
 and
\begin{eqnarray*}
&{}& \int_0^\infty e^{-\rho^2/4t} e^{-\lambda^2/4t}
\int_0^\infty  w_\lambda(v) {dt \o t^{n/2}} \\
&=& \Gamma({n \o 2}-1) 2^{n-2} \int_0^\infty
{ w_\lambda(v) dv \o (\lambda^2+\rho^2)^{{n \o 2}-1}}\/.
\end{eqnarray*}
Taking into account appropriate constants, we obtain the formulas for the Poisson kernel.
\end{proof}

\begin{thm} \label{Pkernel1}
  $$ \lim_{|y|\to\infty} |y|^{n+2\mu -1}P_1(x,y)=c_\mu, \quad \mu>0;$$
   $$ \lim_{|y|\to\infty} (\log|y|)^2|y|^{n -1}P_1(x,y)=c_0, \quad \mu=0;$$
  for some positive $c_\mu$.
   \end{thm}
   \begin{proof}
    From Corollary \ref{Poisson}, Theorem \ref{Ainfty} and arguments presented in its proof
    (boundedeness of $t^{\mu+1} q_\mu(t)$ for $\mu>0$ or $(\log t)^2 q_\mu(t)$,
     for $\mu=0$)
    we obtain that the asymptotic behaviour of $P_1(x,y)$ is the same
    (up to a positive constant) as of
   $$\int_0^\infty t^{-1-\mu-(n-1)/2} e^{-|y|^2/4t} dt = 2^{2\mu+n-1}
   |y|^{-2\mu-n+1} \int_0^\infty u^{\mu+(n-3)/2} e^{-u} du\/,
   $$
   for $\mu>0$. Similar arguments apply for $\mu=0$.
   \end{proof}

  \section*{Acknowledgements}
  The authors would like to thank
  T. Zak for stimulating conversations on the
  subject.

\end{document}